\documentclass[12pt,reqno]{amsart}
\usepackage{amsmath, amsthm, amscd, amsfonts, amssymb, graphicx, color}
\usepackage[bookmarksnumbered, colorlinks, plainpages]{hyperref}

\theoremstyle{definition}

\theoremstyle{remark}

\numberwithin{equation}{section}
\begin{document}
\begin{center}
{\textbf{Geometrical structure in a perfect fluid spacetime with conformal Ricci-Yamabe soliton}}
\end{center}
\vskip 0.3cm
\begin{center}By\end{center}\vskip 0.3cm
\begin{center}
{Soumendu Roy \footnote{The first author is the corresponding author, supported by Swami Vivekananda Merit Cum Means Scholarship, Government of West Bengal, India.}, Santu Dey $^2$ and Arindam~~Bhattacharyya $^3$}
\end{center}
\vskip 0.3cm
\address[Soumendu Roy]{Department of Mathematics,Jadavpur University, Kolkata-700032, India}
\email{soumendu1103mtma@gmail.com}

\address[Santu Dey]{Department of Mathematics, Bidhan Chandra College, Asansol - 4, West Bengal-713304 , India}
\email{santu.mathju@gmail.com}

\address[Arindam Bhattacharyya]{Department of Mathematics,Jadavpur University, Kolkata-700032, India}
\email{bhattachar1968@yahoo.co.in}
\vskip 0.5cm
\begin{center}
\textbf{Abstract}\end{center}
The present paper is to deliberate the geometric composition of a perfect fluid spacetime with torse-forming vector field $\xi$ in connection with conformal Ricci-Yamabe metric and conformal $\eta$-Ricci-Yamabe metric. Here we have delineated the conditions for conformal Ricci-Yamabe soliton to be expanding, steady or shrinking. Later, we have acquired Laplace equation from conformal $\eta$-Ricci-Yamabe soliton equation when the potential vector field $\xi$ of the soliton is of gradient type. Lastly, we have designated perfect fluid with Robertson-Walker spacetime and some applications of physics and gravity.
\vspace {0.3cm}
\section{\textbf{Introduction}}
In 1982, R. S. Hamilton \cite{rsham} introduced the concept of Ricci flow, which is an evolution equation for metrics on a Riemannian manifold. The Ricci flow equation is given by:
\begin{equation}\label{1.1}
\frac{\partial g}{\partial t} = -2S
\end{equation}
on a compact Riemannian manifold $M$ with Riemannian metric $g$.\\
A self-similar solution to the Ricci flow \cite{rsham}, \cite{topping} is called a Ricci soliton \cite{rsha} if it moves only by a one parameter family of diffeomorphism and scaling. The Ricci soliton equation is given by:
\begin{equation}\label{1.2}
\pounds_V g + 2S + 2\Lambda g=0,
\end{equation}
where $\pounds_V$ is the Lie derivative in the direction of $V$, $S$ is Ricci tensor, $g$ is Riemannian metric, $V$ is a vector field and $\Lambda$ is a scalar. The Ricci soliton is said to be shrinking, steady and expanding accordingly as $\Lambda$ is negative, zero and positive respectively.\\
In 2015, N. Basu and A. Bhattacharyya \cite{nbab} established the notion of conformal Ricci soliton \cite{soumendu}, \cite{roy} as:
\begin{equation}\label{1.3}
\pounds_V g + 2S + \Big[2\Lambda - \Big(p + \frac{2}{n}\Big)\Big]g=0,
\end{equation}
where $S$ is the Ricci tensor, $p$ is a scalar non-dynamical field(time dependent scalar field),  $\Lambda$ is constant, $n$ is the dimension of the manifold. \\
The notion of Conformal $\eta$- Ricci soliton was introduced by Mohd Danish Siddiqi \cite{mohd} in 2018, which can be written as:
\begin{equation}\label{1.4}
  \pounds_\xi g + 2S +\Big[2\Lambda - \Big(p + \frac{2}{n}\Big)\Big] g+2 \mu \eta \otimes \eta=0,
\end{equation}
where  $\pounds_\xi$ is the Lie derivative along the vector field $\xi$ , $\Lambda$, $\mu$ are contants, $S$, $p$, $n$ are same as defined in \eqref{1.3}.\\
The concept of Yamabe flow was first introduced by Hamilton \cite{rsha} to construct Yamabe metrics on compact Riemannian manifolds. On a Riemannian or pseudo-Riemannian manifold $M$, a time-dependent metric $g(\cdot, t)$ is said to evolve by the Yamabe flow if the metric $g$ satisfies the given equation,
\begin{equation}\label{1.5}
  \frac{\partial }{\partial t}g(t)=-rg(t),\hspace{0.5cm} g(0)=g_{0},
\end{equation}
where $r$ is the scalar curvature of the manifold $M$.\\\\
In 2-dimension the Yamabe flow is equivalent to the Ricci flow \cite{rsham} (defined by $\frac{\partial }{\partial t}g(t) = -2S(g(t))$, where $S$ denotes the Ricci tensor). But in dimension $> 2$ the Yamabe and Ricci flows do not agree, since the Yamabe flow preserves the conformal class of the metric but the Ricci flow does not in general.\\
A Yamabe soliton \cite{barbosa} correspond to self-similar solution of the Yamabe flow, is defined on a Riemannian or pseudo-Riemannian manifold $(M, g)$ as:
\begin{equation}\label{1.6}
  \frac{1}{2}\pounds_V g = (r-\Lambda)g,
\end{equation}
where $\pounds_V g$ denotes the Lie derivative of the metric $g$ along the vector field $V$, $r$ is the scalar curvature and $\Lambda$ is a constant. Moreover a Yamabe soliton is said to be expanding, steady, shrinking depending on $\Lambda$ being positive, zero, negative respectively. If $\Lambda$ is a smooth function then \eqref{1.6} is called almost Yamabe soliton \cite{barbosa}. \\
Since the introduction of Ricci soliton and Yamabe soliton, many authors ( \cite{roy2}, \cite{roy3}, \cite{ghosh}, \cite{dong}, \cite{joma}, \cite{singh}) have studied these solitons on contact manifolds.\\\\
Recently in 2019, S. G\"uler and M. Crasmareanu \cite{guler} introduced a new geometric flow which is a scalar combination of Ricci and Yamabe flow under the name Ricci-Yamabe map. This flow is also known as Ricci-Yamabe flow of the type $(\alpha, \beta)$.\\
Let $(M^n, g)$ be a Riemannian manifold and $T^{s}_2 (M)$ be the linear space of its symmetric tensor fields of (0, 2)-type and $Riem(M) \subsetneqq T^{s}_2 (M)$  be the infinite space of its Riemannian metrics. In \cite{guler}, the authors have stated the following definition:\\\\
\textbf{Definition 1.1:}\cite{guler} A Riemannian flow on $M$ is a smooth map:
 $$g:I\subseteq \mathbb{R}\rightarrow Riem(M),$$
 where $I$ is a given open interval. We can call it also as time-dependent (or non-stationary) Riemannian metric.\\
\textbf{Definition 1.2:}\cite{guler} The map $RY^{(\alpha, \beta, g)}: I \rightarrow T^{s}_2 (M)$ given by:
$$ RY^{(\alpha, \beta, g)}:= \frac{\partial }{\partial t}g(t)+2\alpha S(t)+\beta r(t)g(t),$$
is called the $(\alpha, \beta)$-Ricci-Yamabe map of the Riemannian flow of the Riemannian flow $(M^n, g)$, where $\alpha, \beta$ are some scalas. If $RY^{(\alpha, \beta, g)} \equiv 0$, then $g(\cdot)$ will be called an $(\alpha, \beta)$-Ricci-Yamabe flow.\\
Also in \cite{guler}, the authors characterized that the $(\alpha, \beta)$-Ricci-Yamabe flow is said to be:\\
$\bullet$ Ricci flow \cite{rsham} if $\alpha = 1$, $\beta = 0$.\\
$\bullet$ Yamabe flow \cite{rsha} if $\alpha = 0$, $\beta = 1$.\\
$\bullet$ Einstein flow \cite{CATINO} if $\alpha = 1$, $\beta = -1$.\\\\
A soliton to the Ricci-Yamabe flow is called Ricci-Yamabe solitons if it moves only by one parameter group of diffeomorphism and scaling. The metric of the Riemannain manifold $(M^n, g)$, $n >2$ is said to admit $(\alpha, \beta)$--Ricci-Yamabe soliton or simply Ricci-Yamabe soliton (RYS) $(g, V, \Lambda,\alpha, \beta)$ if it satisfies the equation:
\begin{equation}\label{1.7}
  \pounds_V g+2\alpha S=[2\Lambda-\beta r]g,
\end{equation}
where $\pounds_V g$ denotes the Lie derivative of the metric $g$ along the vector field $V$, $S$ is the Ricci tensor, $r$ is the scalar curvature and $\Lambda, \alpha,\beta$ are real scalars.\\\\
In the above equation if the vector field $V$ is the gradient of a smooth function $f$ (denoted by $Df$, $D$ denotes the gradient operator) then the equation \eqref{1.7} is called gradient Ricci-Yamabe soliton (GRYS) and it is defined as:
\begin{equation}\label{1.8}
  Hess f+ \alpha S=\Big[\Lambda-\frac{1}{2}\beta r \Big]g,
\end{equation}
where $Hess f$ is the Hessian of the smooth function $f$.\\
Moreover the Ricci-Yamabe soliton and gradient Ricci-Yamabe soliton  is said to be expanding, steady or shrinking according as $\Lambda$ is negative, zero, positive respectively. Also if $\Lambda, \alpha, \beta$ become smooth function then \eqref{1.7} and \eqref{1.8} are called almost Ricci-Yamabe soliton and gradient almost Ricci-Yamabe soliton respectively\\\\
Now using \eqref{1.7} and \eqref{1.3}, we introduce the notion of conformal Ricci-Yamabe soliton as:\\\\
\textbf{Definition 1.3:} A Riemannian manifold $(M^n,g)$, $n>2$ is said to admit conformal Ricci-Yamabe soliton if
\begin{equation}\label{1.9}
  \pounds_V g+2\alpha S+\Big[2\Lambda-\beta r -\Big(p + \frac{2}{n}\Big)\Big]g=0,
\end{equation}
where $\pounds_V g$ denotes the Lie derivative of the metric $g$ along the vector field $V$, $S, r, \Lambda, \alpha,\beta$ are same as defined in \eqref{1.7} and $p, n$ are same as defined in \eqref{1.3}.\\
The conformal Ricci-Yamabe soliton is said to be expanding, steady, shrinking depending on $\Lambda$ being positive, zero, negative respectively. If the vector field $V$ is of gradient type i.e $V = grad(f)$, for $f$ is a smooth function on $M$, then the equation \eqref{1.9} is called conformal gradient Ricci-Yamabe soliton.\\\\
Also using \eqref{1.7} and \eqref{1.4}, we develop the notion of conformal $\eta$-Ricci-Yamabe soliton as:\\\\
\textbf{Definition 1.4:} A Riemannian manifold $(M^n,g)$, $n>2$ is said to admit conformal $\eta$-Ricci-Yamabe soliton if
\begin{equation}\label{1.10}
  \pounds_\xi g+2\alpha S+\Big[2\Lambda-\beta r -\Big(p + \frac{2}{n}\Big)\Big]g+2 \mu \eta \otimes \eta=0,
\end{equation}
where $\pounds_\xi g$ denotes the Lie derivative of the metric $g$ along the vector field $\xi$, $\Lambda, \mu$ are contants, $S, r, \alpha,\beta$ are same as defined in \eqref{1.7} and $p, n$ are same as defined in \eqref{1.3}.\\
If the vector field $\xi$ is of gradient type i.e $\xi = grad(f)$, for $f$ is a smooth function on $M$, then the equation \eqref{1.10} is called conformal gradient $\eta$-Ricci-Yamabe soliton.\\\\
On the other side, in 1915, Albert Einstein established General relativity, also known as the general theory of relativity (GTR), which is the geometric theory of gravitation. In this theory, the gravitational field is the spacetime curvature and its source is energy–momentum tensor. The goal to develop differential geometry and relativistic fluids model in the mathematical language are most efficient for understanding general relativity. The spacetime of general relativity and cosmology can be modeled as a connected 4-dimensional Lorentzian manifold which is a special subclass of pseudo-Riemannian manifolds with Lorentzian metric $g$ with signature (-, +, +, +) has great importance in general relativity. The geometry of Lorentzian manifold begins with the study of causal character of vectors of the manifold, due to this causality that Lorentzian manifold becomes a convenient choice for the study of general relativity.\\
The energy–momentum tensor plays the importent role as a matter content of the spacetime, matter is assumed to be fluid having density, pressure and having dynamical and kinematic quantities like velocity, acceleration, vorticity, shear and expansion \cite{ahsan}, \cite{stephani}. The matter content of the universe is assumed to act like a perfect fluid in standard cosmological models. The most suitable example of perfect fluid is dust fluid.\\\\\\
The outline of the article goes as follows:\\
In section 2, after a brief introduction, we have discussed some needful properties of perfect fluid which will be used in the later sections. Section 3 deals with some applications of conformal Ricci-Yamabe soliton structure in perfect fluid spacetime with torse-forming vector field. In this section we have contrived the conformal Yamabe soliton in perfect fluid spacetime with torse forming vector field to accessorize the nature of this soliton on the mentioned spacetime. We have also considered the potential vector field $V$ of the solition as conformal Killing vector field to characterized the vector field. Section 4 is devoted to form the Laplace equation from conformal $\eta$-Ricci-Yamabe soliton equation when the potential vector field $\xi$ of the soliton is of gradient type. In section 5 and section 6, we have shown the physical connection of perfect fluid with Robertson-Walker spacetime and the application of Laplace equation in physics and gravity respectively.\\
\vspace {0.3cm}
\section{\textbf{Perfect fluid spacetime with torse-forming vector field}}
A perfect fluid is a fluid that can be completely characterized by its rest frame mass density and isotropic pressure. A perfect fluid has no shear stress, viscosity or heat conduction and it is distinguished by an energy-momentum tensor $T$ of the form \cite{neill}:
\begin{equation}\label{2.1}
  T(X,Y)=\rho g(X,Y)+(\sigma+\rho)\eta(X)\eta(Y),
\end{equation}
where $\rho, \sigma$ are the isotropic pressure and energy-density respectively and $\eta(X) = g(X,\xi)$ is 1-form, which is equivalent to the unit vector $\xi$ and $g(\xi,\xi) = -1$. The field equation governing the perfect fluid motion is Einstein's gravitational equation \cite{neill}:
\begin{equation}\label{2.2}
  S(X,Y)+\Big[\lambda-\frac{r}{2} \Big]g(X,Y)=\kappa T(X,Y),
\end{equation}
where $S, r$ are the Ricci tensor and scalar curvature of $g$ respectively, $\lambda$ is the cosmological constant and $\kappa$ is the gravitational constant, which can be considered as $8\pi G$, where $G$ is the universal gravitational constant.\\
Using \eqref{2.1}, the above equation takes the form:
\begin{equation}\label{2.3}
  S(X,Y)=\Big[-\lambda+\frac{r}{2}+\kappa \rho \Big]g(X,Y)+\kappa(\sigma+\rho)\eta(X)\eta(Y).
\end{equation}
Let $(M^4, g)$ be a relativistic perfect fluid spacetime which satisfies \eqref{2.3}. Then by contracting \eqref{2.3} and considering $g(\xi, \xi) = -1$, we obtain,
\begin{equation}\label{2.4}
r=4\lambda+\kappa(\sigma-3\rho).
\end{equation}
Using the value of $r$ from the above equation, \eqref{2.3} becomes,
\begin{equation}\label{2.5}
  S(X,Y)=\Big[\lambda+\frac{\kappa(\sigma-\rho)}{2} \Big]g(X,Y)+\kappa(\sigma+\rho)\eta(X)\eta(Y).
\end{equation}
Hence the Ricci operator $Q$ can be written as:
\begin{equation}\label{2.6}
  QX= \Big[\lambda+\frac{\kappa(\sigma-\rho)}{2} \Big]X+\kappa(\sigma+\rho)\eta(X)\xi,
\end{equation}
where $g(QX, Y)= S(X, Y)$.\\\\
\textbf{Example 2.1:} A radiation fluid is a perfect fluid with $\sigma = 3\rho$ and so the energy momentum tensor $T$ becomes,
\begin{equation}\label{2.6new}
   T(X,Y)=\rho [g(X,Y)+4 \eta(X)\eta(Y)],
\end{equation}
From \eqref{2.4}, we can say that a radiation fluid has constant scalar curvature $r$ equal to $4\lambda$.\\\\
Now we take a special case when $\xi$ is a torse-forming vector field \cite{blaga}, \cite{yano} of the form:
\begin{equation}\label{2.7}
  \nabla_X \xi = X+\eta(X)\xi.
\end{equation}
Also on a perfect fluid spacetime if the vector field $\xi$ is torse-forming, then the following relations hold \cite{blaga}:
\begin{equation}\label{2.8}
  \nabla_\xi \xi =0,
\end{equation}
\begin{equation}\label{2.9}
  (\nabla_X \eta)(Y)=g(X,Y)+\eta(X)\eta(Y)
\end{equation}
\begin{equation}\label{2.10}
  R(X,Y)\xi=\eta(Y)X-\eta(X)Y,
\end{equation}
\begin{equation}\label{2.11}
  \eta(R(X,Y)Z)=\eta(X)g(Y,Z)-\eta(Y)g(X,Z),
\end{equation}
for all vector fields $X, Y, Z$.\\
Using \eqref{2.7}, we have,
\begin{eqnarray}\label{2.12}
  (\pounds_\xi g)(X,Y) &=& g(\nabla_X \xi,Y)+g(X,\nabla_Y \xi) \nonumber\\
                        &=& 2 [g(X,Y)+\eta(X)\eta(Y)],
\end{eqnarray}
for all vector fields $X, Y$.\\\\
Perfect fluid are frequently used in general relativity to model idealized distribution of matter such as the interior of a star or an isotropic universe. In general relativity and symmetries of spacetime one often employs a perfect fluid energy momentum tensor \eqref{2.1} to represent the source of the gravitational field. A perfect fluid has two thermodynamic degrees of freedom.\\\\
In general relativity, a perfect fluid solution is an exact solution of the Einstein field equation in which the gravitational field is produced entirely by the mass., momentum and stress density of a fluid.\\
In astrophysics, fluid solutions are often employed as stellar models. It might help to think of a perfect gas as a special case of a perfect fluid. In cosmology, fluid solutions are often used as cosmological models.\\
There are some special case of fluid solutions:\\\\
(i) A dust is a pressureless perfect fluid with the energy momentum tensor $T(X,Y)=\sigma \eta(X)\eta(Y)$.\\\\
(ii) A radiation fluid is a perfect fluid with \eqref{2.6new}.\\\\
These two are often used as cosmological models for matter dominated and radiation dominated epochs. While in general it requires ten functions to specify a fluid,  a perfect fluid requires only two whereas dust and radiation fluid each requires only one function.\\
It is much easier to find such solutions than it is to find a general fluid solution.\\\\
Among the perfect fluids other than dust or radiation fluids, by far the most important special case is that of the static spherically symmetric perfect fluid solutions. These can always be matched to a schwarzschild vaccum across a spherical surfaces, so they can be used as interior solutions in a stellar model. Also the characteristic polynomial of the Einstein tensor in a perfect fluid must have the form:
$$\chi(\tau)=(\tau-8\pi\sigma)(\tau-8\pi\rho)^3,$$
 where $\sigma, \rho$ are the density and pressure of fluid respectively.\\
Perfect fluid solutions which feature positive pressure include various radiation fluid models from cosmology, including\\\\
(a) FRW radiation fluids, often referred to us the radiation dominated FRW models.\\\\
(b) Wahlquist fluid, which has similar symmetries to the kerr vacuum, leading to initial hopes that it might provide the interior solutions for a simple model of a rotating star.\\\\
(c) The equation of state of the perfect fluid may be used in Friedmann- Lemaitre- Robertson Walker equations to describe the evolution of the universe.\\\\
\section{\textbf{Conformal Ricci-Yamabe soliton structure in perfect fluid spacetime with torse-forming vector field}}
In this section, we study conformal Ricci-Yamabe soliton structure in a perfect fluid spacetime whose timelike velocity vector field $\xi$ is torse-forming.\\
Taking $V$ as a torse-forming vector field $\xi$ in the soliton equation \eqref{1.9} and putting $n = 4$, we get,
\begin{equation}\label{3.1}
  (\pounds_\xi g)(X,Y)+2\alpha S(X,Y)+\Big[2\Lambda-\beta r -\Big(p + \frac{1}{2}\Big)\Big]]g(X,Y)=0.
\end{equation}
Using \eqref{2.12}, the above equation becomes,
\begin{equation}\label{3.2}
  2 [g(X,Y)+\eta(X)\eta(Y)]+2\alpha S(X,Y)+\Big[2\Lambda-\beta r -\Big(p + \frac{1}{2}\Big)\Big]]g(X,Y)=0.
\end{equation}
In view of \eqref{2.5}, we obtain,
\begin{equation}\label{3.3}
  \Big[\Lambda-\frac{\beta r}{2}-\frac{1}{2}\Big(p + \frac{1}{2}\Big)+\alpha\lambda+\frac{\alpha\kappa(\sigma-\rho)}{2}+1\Big]+\Big[\alpha\kappa(\sigma+\rho)+1\Big]\eta(X)\eta(Y)=0.
\end{equation}
Taking $X= Y= \xi$ in the above equation, we get,
\begin{equation}\label{3.4}
  \Lambda=\frac{\alpha\kappa(\sigma+3\rho)}{2}+\frac{\beta r}{2}-\alpha\lambda+\frac{1}{2}\Big(p + \frac{1}{2}\Big).
\end{equation}
Using \eqref{2.4}, we have,
\begin{equation}\label{3.5}
  \Lambda=\frac{\kappa}{2}\Big[(\alpha+\beta)\sigma+3(\alpha-\beta)\rho \Big]+(2\beta-\alpha)\lambda+\frac{1}{2}\Big(p + \frac{1}{2}\Big).
\end{equation}
So we can state the following:\\\\
\textbf{Theorem 3.1.} {\em If a perfect fluid spacetime with torse-forming vector field $\xi$ admits a conformal Ricci-Yamabe soliton $(g, \xi, \Lambda,\alpha, \beta)$, then the soliton is expanding, steady, shrinking according as, $\frac{\kappa}{2}\Big[(\alpha+\beta)\sigma+3(\alpha-\beta)\rho \Big]+(2\beta-\alpha)\lambda+\frac{1}{2}\Big(p + \frac{1}{2}\Big) \gtreqqless 0$. }\\\\
\textbf{Remark 3.2:} {\em In \eqref{3.5}, if we take $\Big(p + \frac{1}{2}\Big)=0$, then $\Lambda=\frac{\kappa}{2}\Big[(\alpha+\beta)\sigma+3(\alpha-\beta)\rho \Big]+(2\beta-\alpha)\lambda$ and in this case the conformal Ricci-Yamabe soliton becomes Ricci-Yamabe soliton and we obtain that the soliton is expanding, steady, shrinking according as,
$\frac{\kappa}{2}\Big[(\alpha+\beta)\sigma+3(\alpha-\beta)\rho \Big]+(2\beta-\alpha)\lambda \gtreqqless 0$. }\\\\
A spacetime symmetry of physical interest is the conformal Killing vector as it preserves the metric up to a conformal factor. A vector field $V$ is said to be a conformal Killing vector field iff the following relation holds:
\begin{equation}\label{3.6}
  (\pounds_V g)(X,Y)=2\Phi g(X,Y),
\end{equation}
where $\Phi$ is some function of the co-ordinates(conformal scalar).\\
Moreover if $\Phi$ is not constant the conformal Killing vector field $V$ is said to be proper. Also when $\Phi$ is constant, $V$ is called homothetic vector field and when the constant $\Phi$ becomes non zero, $V$ is said to be proper homothetic vector field. If $\Phi = 0$ in the above equation, then $V$ is called Killing vector field.\\
Let us assume that in the equation \eqref{1.9}, the potential vector field $V$ is conformal Killing vector field. Then using \eqref{3.6} and \eqref{1.9}, we get,
\begin{equation}\label{3.7}
  \alpha S(X,Y)=-\Big[\Lambda+\Phi-\frac{\beta r}{2}-\frac{1}{2}\Big(p + \frac{1}{2}\Big)\Big]g(X,Y),
\end{equation}
which leads to the fact that the spacetime is Einstein, provided $\alpha \neq 0$.\\
Conversely, assuming that the perfect fluid spacetime with torse-forming vector filed $\xi$ is Einstein spacetime,i.e. $S(X,Y) = \theta g(X,Y)$.\\
Then the equation \eqref{1.9} becomes,
\begin{equation}\label{3.8}
  (\pounds_V g)(X,Y)=-\Big[2\Lambda+2\alpha \theta-\beta r -\Big(p + \frac{1}{2}\Big)\Big]]g(X,Y),
\end{equation}
which can be written as,
\begin{equation}\label{3.9}
  (\pounds_V g)(X,Y)=2\Psi g(X,Y),
\end{equation}
where $\Psi = -\Big[\Lambda+\alpha \theta-\frac{\beta r}{2} -\frac{1}{2}\Big(p + \frac{1}{2}\Big)\Big]$.\\
Thus from \eqref{3.9}, $V$ becomes a conformal Killing vector field.\\
Hence we can state the following:\\\\
\textbf{Theorem 3.3} {\em Let a perfect fluid spacetime with torse-forming vector field $\xi$ admit a conformal Ricci-Yamabe soliton $(g, V, \Lambda,\alpha, \beta)$. The potential vector field $V$ is a conformal Killing vector field iff the spacetime is Einstein, provided $\alpha \neq 0$.}\\\\
Now in view of \eqref{3.7} and \eqref{2.5}, we obtain,
\begin{equation}\label{3.10}
  \Big[\Lambda+\Phi+\alpha\lambda+\frac{\alpha\kappa(\sigma-\rho)}{2}-\frac{\beta r}{2} -\frac{1}{2}\Big(p + \frac{1}{2}\Big)\Big]g(X,Y)+\Big[\alpha\kappa(\sigma+\rho)\Big]\eta(X)\eta(Y)=0.
\end{equation}
Taking $Y = \xi$ in the above equation and considering $\eta(\xi) = -1$, we have,
\begin{equation}\label{3.11}
  \Big[\Lambda+\Phi+\alpha\lambda-\frac{\alpha\kappa(\sigma+3\rho)}{2}-\frac{\beta r}{2} -\frac{1}{2}\Big(p + \frac{1}{2}\Big)\Big]\eta(X)=0.
\end{equation}
Since $\eta(X) \neq 0$, then we get,
\begin{equation}\label{3.12}
  \Lambda+\Phi+\alpha\lambda-\frac{\alpha\kappa(\sigma+3\rho)}{2}-\frac{\beta r}{2} -\frac{1}{2}\Big(p + \frac{1}{2}\Big)=0.
\end{equation}
Substituting the value of $r$ from \eqref{2.4}, the above equation reduces to,
\begin{equation}\label{3.13}
  \Phi=\frac{\kappa}{2}\Big[(\alpha+\beta)\sigma+3(\alpha-\beta)\rho \Big]+(2\beta-\alpha)\lambda-\Lambda+\frac{1}{2}\Big(p + \frac{1}{2}\Big).
\end{equation}
Hence we can state the following:\\\\
\textbf{Theorem 3.4.} {\em Let a perfect fluid spacetime with torse-forming vector field $\xi$ admit a conformal Ricci-Yamabe soliton $(g, V, \Lambda,\alpha, \beta)$. The potential vector field $V$ is a conformal Killing vector field, then $V$ is\\
(i) proper conformal Killing vector field if $\alpha, \beta, p$ are not constant.\\
(ii) homothetic vector field if $\alpha,\beta, p$ are constant.}\\\\
Using the property of Lie derivative we can write,
\begin{equation}\label{3.14}
  (\pounds_V g)(X,Y)=g(\nabla_X V, Y)+g(\nabla_Y V, X),
\end{equation}
for any vector fields $X, Y$.\\
Then from \eqref{2.5} and \eqref{3.14}, \eqref{1.9} takes the form,
\begin{multline}\label{3.15}
   g(\nabla_X V, Y)+g(\nabla_Y V, X)+\Big[2\Lambda-\beta r -\Big(p + \frac{1}{2}\Big)+2\alpha\Big\{\lambda+\frac{\kappa(\sigma-\rho)}{2}\Big\}\Big]g(X,Y) \\
  +2\alpha\kappa(\sigma+\rho)\eta(X)\eta(Y)=0.
\end{multline}
Suppose $\omega$ is a 1-form, which is metrically equivalent to $V$ and is given by $\omega(X) = g(X,V)$ for an arbitrary vector field $X$. Then the exterior derivative $d\omega$ of $\omega$ can be written as:
\begin{equation}\label{3.16}
  2(d\omega)(X,Y)=g(\nabla_X V,Y)-g(\nabla_Y V,X).
\end{equation}
As $d\omega$ is skew-symmetric, so if we define a tensor field $F$ of type (1,1) by,
\begin{equation}\label{3.17}
  (d\omega)(X,Y)=g(X,FY),
\end{equation}
then $F$ is skew self-adjoint i.e. $g(X, FY)=-g(FX, Y)$.\\
So \eqref{3.17} can be written as:
\begin{equation}\label{3.18}
  (d\omega)(X,Y)=-g(FX,Y)
\end{equation}
Using \eqref{3.18}, \eqref{3.16} becomes,
\begin{equation}\label{3.19}
  g(\nabla_X V,Y)-g(\nabla_Y V,X)=-2g(FX,Y).
\end{equation}
Adding \eqref{3.19} and \eqref{3.15} side by side and factoring out $Y$, we get,
\begin{equation}\label{3.20}
  \nabla_X V=-FX-\Big[\Lambda-\frac{\beta r}{2} -\frac{1}{2}\Big(p + \frac{1}{2}\Big)+\alpha\Big\{\lambda+\frac{\kappa(\sigma-\rho)}{2}\Big\}\Big]X-\alpha\kappa(\sigma+\rho)\eta(X)\xi.
\end{equation}\\
Substituting the above equation in $R(X,Y)V=\nabla_X \nabla_Y V - \nabla_Y \nabla_X V - \nabla_{[X,Y]} V$, we have,
\begin{equation}\label{3.21}
  R(X,Y)V=(\nabla_Y F)X-(\nabla_X F)Y+\alpha\kappa(\sigma+\rho)[Y\eta(X)-X\eta(Y)].
\end{equation}
Noting that $d\omega$ is closed, we obtain,
\begin{equation}\label{3.22}
  g(X,(\nabla_Z F)Y)+g(Y,(\nabla_X F)Z)+g(Z,(\nabla_Y F)X)=0.
\end{equation}
Making inner product of \eqref{3.21} with respect to $Z$, we get,
\begin{eqnarray}\label{3.23}
  g(R(X,Y)V,Z) &=& g((\nabla_Y F)X,Z)-g((\nabla_X F)Y,Z) \nonumber \\
               &+& \alpha\kappa(\sigma+\rho)[g(Y,Z)\eta(X)-g(X,Z)\eta(Y)].
\end{eqnarray}
As $F$ is skew self-adjiont, then $\nabla_X F$ is also skew self-adjiont. Then using \eqref{3.22}, \eqref{3.23} takes the form,
\begin{equation}\label{3.24}
  g(R(X,Y)V,Z)=\alpha\kappa(\sigma+\rho)[g(Y,Z)\eta(X)-g(X,Z)\eta(Y)]-g(X,(\nabla_Z F)Y).
\end{equation}
Putting $X = Z =e_i$ in the above equation, where $e_i$'s are a local orthonormal frame and summing over $i=1,2,3,4$, we obtain,
\begin{equation}\label{3.25}
  S(Y,V)=-3\alpha\kappa(\sigma+\rho)\eta(Y)-(div F)Y,
\end{equation}
where $div F$ is the divergence of the tensor field $F$.\\
Equating \eqref{2.5} and \eqref{3.25}, we get,
\begin{equation}\label{3.26}
  (div F)Y=-\kappa(\sigma+\rho)[3\alpha+\eta(V)]\eta(Y)-\Big[\lambda+\frac{\kappa(\sigma-\rho)}{2} \Big]\omega(Y).
  \end{equation}
Now we compute the covariant derivative of the squared $g$-norm of $V$ using \eqref{3.20} as follows:
\begin{eqnarray}\label{3.27}
  \nabla_X \mid V\mid^2&=&2g(\nabla_X V,V) \nonumber\\
  &=&-2g(FX,V)-\Big[2\Lambda-\beta r -\Big(p + \frac{1}{2}\Big) \nonumber\\
  &+&2\alpha\Big\{\lambda+\frac{\kappa(\sigma-\rho)}{2}\Big\}\Big]g(X,V)- 2\alpha\kappa(\sigma+\rho)\eta(X)\eta(V).
\end{eqnarray}
From \eqref{2.5}, \eqref{1.9} becomes,
\begin{eqnarray}\label{3.28}
 (\pounds_V g)(X,Y) &=& - \Big[2\Lambda-\beta r -\Big(p + \frac{1}{2}\Big)+2\alpha\Big\{\lambda+\frac{\kappa(\sigma-\rho)}{2}\Big\}\Big]g(X,Y)\nonumber\\
  &-& 2\alpha\kappa(\sigma+\rho)\eta(X)\eta(Y).
\end{eqnarray}
Using the above equation, \eqref{3.27} takes the form,
\begin{equation}\label{3.29}
  \nabla_X \mid V\mid^2+2g(FX,V)-(\pounds_V g)(X,V)=0.
\end{equation}
So we can state the following:\\\\
\textbf{Theorem 3.5} {\em If a perfect fluid spacetime with torse-forming vector field $\xi$ admits a conformal Ricci-Yamabe soliton $(g, V, \Lambda,\alpha, \beta)$, then the vector $V$ and its metric dual 1-form $\omega$ satisfies the relation $$(div F)Y=-\kappa(\sigma+\rho)[3\alpha+\eta(V)]\eta(Y)-\Big[\lambda+\frac{\kappa(\sigma-\rho)}{2} \Big]\omega(Y)$$ and $$\nabla_X \mid V\mid^2+2g(FX,V)-(\pounds_V g)(X,V)=0.$$}\\\\
\section{\textbf{Conformal $\eta$-Ricci-Yamabe soliton structure in perfect fluid spacetime}}
Let $(M^4,g)$ be a general relativistic perfect fluid spacetime and $(g, \xi, \Lambda,\mu,\alpha, \beta)$ be a conformal $\eta$-Ricci-Yamabe soliton in $M$.\\
Then writting explicitly the Lie derivative $(\pounds_\xi g)$ as $ (\pounds_\xi g)(X,Y)= g(\nabla_X \xi, Y) + g(X, \nabla_Y \xi)$ and from \eqref{1.10} and \eqref{2.5}, we obtain,
\begin{multline}\label{4.1}
  g(\nabla_X \xi, Y) + g(X, \nabla_Y \xi)+2\alpha\Big[\Big\{\lambda+\frac{\kappa(\sigma-\rho)}{2} \Big\}g(X,Y)+\kappa(\sigma+\rho)\eta(X)\eta(Y)\Big] \\
 +\Big[2\Lambda-\beta r -\Big(p + \frac{1}{2}\Big)\Big]g(X,Y)+2 \mu \eta(X)\eta(Y)=0,
\end{multline}
for any vector fields $X, Y$.\\
Then the above equation can be written as,
\begin{multline}\label{4.2}
 \Big[\Lambda-\frac{\beta r}{2}-\frac{1}{2}\Big(p + \frac{1}{2}\Big)+\alpha\lambda+\frac{\alpha\kappa(\sigma-\rho)}{2}\Big]g(X,Y)+\Big[\mu+\alpha\kappa(\sigma+\rho)\Big]\eta(X)\eta(Y) \\
 +\frac{1}{2}\Big
  [g(\nabla_X \xi, Y) + g(X, \nabla_Y \xi)\Big]=0.
\end{multline}
Consider $\{e_i\}_{1 \leq i\leq 4}$ an orthonormal frame field and $\xi = \sum_{i=1}^{4} \xi^i e_i$. We have from \cite{blaga}, $\sum_{i=1}^{4} \epsilon_{ii} (\xi^i)^2=-1$ and $\eta(e_i)=\epsilon_{ii}\xi^i$.\\
Multiplying \eqref{4.2} by $\epsilon_{ii}$ and summing over $i$ for $X = Y = e_i$, we obtain,
\begin{equation}\label{4.3}
  4\Lambda-\mu=4(2\beta-\alpha)\lambda+\kappa(2\beta-\alpha)(\sigma-3\rho)+2\Big(p + \frac{1}{2}\Big)-div(\xi),
\end{equation}
where $div(\xi)$ is the divergence of the vector field $\xi$.\\
Putting $X = Y= \xi$ in \eqref{4.2}, we get,
\begin{equation}\label{4.4}
  \Lambda-\mu=(2\beta-\alpha)\lambda+\frac{\kappa}{2}\Big[(2\beta+\alpha)\sigma-3(2\beta-\alpha)\rho\Big]+\frac{1}{2}\Big(p + \frac{1}{2}\Big).
\end{equation}
Then calculating $\Lambda, \mu$ from \eqref{4.3} and \eqref{4.4}, we get,
\begin{equation}\label{4.5}
  \Lambda=(2\beta-\alpha)\lambda+\frac{\kappa}{2}\Big[\Big(\frac{2\beta-3\alpha}{3}\Big)\sigma-(2\beta-\alpha)\rho\Big]+\frac{1}{2}\Big(p + \frac{1}{2}\Big)-\frac{div(\xi)}{3},
\end{equation}
and
\begin{equation}\label{4.6}
  \mu=-\kappa\Big[\Big(\frac{2\beta+3\alpha}{3}\Big)\sigma-(2\beta-\alpha)\rho\Big]-\frac{div(\xi)}{3}.
\end{equation}
Then we can state the following:\\\\
\textbf{Theorem 4.1} {\em Let $(M^4,g)$ be a $4$-dimensional pseudo-Riemannian manifold and $\eta$ be the $g$-dual 1-form of the gradient vector field $\xi := grad(f)$, with $g(\xi, \xi)=-1$, where $f$ is a smooth function. If $(g, \xi, \Lambda,\mu,\alpha, \beta)$ is a conformal $\eta$-Ricci-Yamabe soliton on $M$, then the Laplacian equation satisfied by $f$ becomes:
\begin{equation}\label{4.7}
  \Delta(f)=-3\Big[\mu+\kappa\Big\{\Big(\frac{2\beta+3\alpha}{3}\Big)\sigma-(2\beta-\alpha)\rho\Big\}\Big].
\end{equation}}\\\\
\textbf{Example 4.2} {\em A conformal $\eta$-Ricci-Yamabe soliton $(g, \xi, \Lambda,\mu,\alpha, \beta)$ in a radiation fluid is given by:
$$\Lambda = (2\beta-\alpha)\lambda-\kappa\alpha\rho+\frac{1}{2}\Big(p + \frac{1}{2}\Big)-\frac{div(\xi)}{3},$$
and
$$\mu = -4\kappa\alpha\rho-\frac{div(\xi)}{3}.$$}\\\\
\section{\textbf{Perfect fluid and Robertson-Walker spacetime}}
Generalized Robertson-Walker(GRW) spacetimes are a natural and wide extension of RW spacetime, where large scale cosmology is staged. They are Lorentzian manifold of dimension $n$ characterized by the metric \cite{mantica},
$$ds^2=-dt^2+q(t)^2g^*_{\gamma \delta}(x_2,x_3,....,x_n)dx^\gamma dx^\delta,\quad\quad \gamma,\delta=2,3,....n,$$
where $g^*_{\gamma \delta}(x_2,x_3,....,x_n)$ is the metric tensor of the Riemannian submanifold and it is the wrapped product $(-I)\times_{q^2} M^*$ (\cite{mantica1}, \cite{alias}), where $(M^*,g^*)$ is a $(n-1)$-dimension Riemannian manifold , $I$ is an interval of the real line and $q>0$ is the smmoth mapping or scale function. If $M^*$ has dimension 3 and has constant curvature, then the spacetime is a Robertson-Walker(RW) spacetime.\\\\
Mantica and Molinari \cite{mantica} have proved that a Lorentzian manifold of dimension $n$ is a GRW spacetime iff it admits a time like torse-forming vector field. If a Lorentzian manifold admits a globally time like vector field, it is called time oriented Lorentzian manifold, physically known as spacetime. Thus the spacetime is a 4-dimensional time oriented Lorentzian manifold.  A Lorentzian manifold has many applications in applied physics, especially in the theory of relativity and cosmology. To study the Lorentzian manifold, the causal character of the vector fields plays an important role and thus it becomes the advantageous choice for the researchers to study the theory of relativity and cosmology.\\\\
Lorentzian manifolds with a Ricci tensor of the form,
$$R_{ij} = Ag_{ij} + Bu_iu_j,$$
where $A$ and $B$ are scalar fields and $u_iu^i = -1$, are often named perfect fluid spacetimes. It is well known that any Robertson-Walker spacetime is a perfect fluid spacetime \cite{neill}, and for $n = 4$, a GRW spacetime is a perfect fluid iff it is a Robertson-Walker spacetime.
So we can establish the fact that, Theorem 3.1, Theorem 3.3, Theorem 3.4, Theorem 3.5 and Theorem 4.1 are also verified on a $4$-dimensional GRW spacetime iff the fluid spacetime is a Robertson-Walker spacetime.\\
The form of the above equation of the Ricci tensor is implied by Einstein’s equation if the energy matter content of space-time is a perfect fluid with velocity vector field $u$. The scalars A and B are linearly related to the pressure and the energy density measured in the locally comoving inertial frame. They are not independent because of the Bianchi identity $\nabla^m R_{im} = \frac{1}{2}\nabla_i R$.\\
Shepley and Taub \cite{shepley} studied a perfect fluid spacetime in dimension $n = 4$, with equation of state and the additional condition that the Weyl tensor has null divergence.\\\\
\section{\textbf{Application of Laplace equation in physics and gravity}}
Laplace equation, a second order P.D.E widely useful in physics as its solution, which is known as harmonic functions occur in problems of electrical, magnetic and gravitational potentials of steady state temparatures and of hydrodynamics.\\\\
$\bullet$ The real and imaginary parts of complex analytic function both satisfy Laplace equation. That is if $z = x+iy$ and $f(x,y)=u(x,y)+iv(x,y)$, then the necessary condition of $f(z)$ to be analytic is that $u$ and $v$ and that be C.R equation be satisfied, $u_x = v_y$, $u_y = -v_x$, where $u_x, u_y$ is the first partial derivatives of $u$ with respect to $x, y$ respectively and $v_x, v_y$ is the first partial derivatives of $v$ with respect to $x, y$ respectively. It follows that $u_{yy}=(-v_x)_y=-(v_y)_x=-(u_{xx}).$\\
Therefore, $u$ satisfies Laplace equation.\\\\
$\bullet$ If we have a region where the charge density is zero (there may be non-zero charge densities at the boundaries), the electric potential $V$ satisfies Laplace equation inside the region. Solving Laplace equation , we get electric potential, which is very important quantity as we can use it to compute the electric field very easily, $E = \nabla V$ and therefore the force $\bar{F} = qE$. There are many interesting cases in physics, where we are concerned with the potential in regions with zero charged density. Classic examples include the region inside and outside a hollow charged sphere , or the region outside charged metal plates. Each of the cases come with different set of boundary conditions on what makes Laplace equation interesting.\\
In general, for a given charged density, $L(x,y,z)$, electric (and gravitational) potentials satisfy poisson's equation, $\nabla^2 V= L(x,y,z)$. Laplace equation or poisson's equation are the simplest examples of a class of P.D.Es called eliptical P.D.Es. A lot of interesting mathematical techniques used to solve electrical P.D.Es are first introduced by Laplace equation.\\\\
$\bullet$ In electrostatics, according to Maxwell's equation, a electric fluid $(u,v)$ in two space dimensions, that is independent of time satisfies,
$$\nabla \times (u,v,0)=(v_x-u_y)\hat{k}=0,$$
and $$\nabla \cdot (u,v) = L,$$ where $L$ is the charge density.\\
The Laplace equation can be used in three dimension problems in electrostatics and fluid flow just as in two dimensions.\\\\
$\bullet$ It has applications in gravity also. Let $\tilde{g}, \tilde{\rho}, G$ be the gravitational field, mass density and gravitational constant. Then Gauss's law for gravitation in differential form is:
$$ \nabla \cdot \tilde{g}=-4\pi G\tilde{\rho}.$$
Also we have, $\nabla^2 V = 4\pi G \tilde{\rho}$, which is poission's equation for gravitational fields.\\
This physical significance is directly equivalent to Theorem 4.1 and \eqref{4.7}, which is a Laplace equation with potential vector field of gradient type.\\
In empty space $\tilde{\rho} = 0$, we have $\nabla^2 V = 0$, which is Laplace equation for gravitational fields.\\\\
 

\begin{thebibliography}{}
\bibitem{alias} L. Alias, A. Romero, M. Sanchez, {\em Compact spacelike hypersurfaces of constant mean curvature in generalized Robertson-Walker spacetimes}, F. Dillen, Geometry and Topology of Submanifolds VII(1995). River Edge NJ, USA: World Scientific, , pp-67-70.
 C. A. Mantica and L. G. Molinari, Generalized Robertson-Walker spacetimes - A survey,
Int. J. Geom. Methods Mod. Phys., 14(2017), 1730001 (27 pages)
\bibitem{ahsan} Z. Ahsan, {\em Tensors: Mathematics of Differential Geometry and Relativity}, PHI Learning Pvt. Ltd, Delhi (2017).
  \bibitem{barbosa} E. Barbosa and E. Ribeiro Jr., {\em On conformal solutions of the Yamabe flow}, Arch. Math.(2013), Vol. 101, pp-79–89.
\bibitem{nbab}Nirabhra Basu and Arindam Bhattacharyya, {\em Conformal Ricci soliton in Kenmotsu manifold}, Global
Journal of Advanced Research on Classical and Modern Geometries(2015), Vol. 4, Isu. 1, pp. 15-21.
\bibitem{blaga} A. M. Blaga, {\em Solitons and geometrical structures in a perfect fluid spacetime}, 	arXiv:1705.04094 [math.DG](2017).
\bibitem{dong} Huai-Dong Cao, Xiaofeng Sun and Yingying Zhang, {\em On the structure of gradient Yamabe solitons}, arXiv:1108.6316v2 [math.DG] (2011).
\bibitem{CATINO} G. Catino and L. Mazzieri, {\em Gradient Einstein solitons}, Nonlinear Anal(2016). Vol. 132, pp-66–94.
\bibitem{joma}Jong Taek Cho, Makoto Kimura,{\em Ricci solitons and real hypersurfaces in a complex space form}, Tohoku Mathematical Journal, Second Series(2009), Vol. 61,Isu. 2,pp. 205-212.
\bibitem{ghosh} Amalendu Ghosh, {\em Yamabe soliton and Quasi Yamabe soliton on Kenmotsu manifold}, Mathematica Slovaca(2020), Vol.70(1), pp-151-160.
\bibitem{guler}  S. G\"uler and M. Crasmareanu, {\em Ricci-Yamabe maps for Riemannian flows and their volume variation and volume entropy}, Turk. J. Math.(2019), Vol.43, pp. 2361-2641.
\bibitem{rsham}R. S. Hamilton, {\em Three Manifold with positive Ricci curvature}, J.Differential Geom.(1982), Vol. 17, Isu.2, pp. 255-306.
\bibitem{rsha} R. S. Hamilton, {\em The Ricci flow on surfaces}, Contemporary Mathematics(1988), Vol. 71, pp. 237-261.
\bibitem{mantica} C. Mantica, L. Molinari and U. C. De, {\em A condition for a perfect-fluid space-time to be a generalized Robertson-Walker space-time}, 	arXiv:1508.05883 [math.DG](2016).
\bibitem{mantica1} C. A. Mantica and L. G. Molinari, {\em Generalized Robertson-Walker spacetimes - A survey},
Int. J. Geom. Methods Mod. Phys., 14(2017), 1730001 (27 pages)
\bibitem{neill} B. O'Neill, {\em Semi-Riemannian Geometry with apllications to Relativity}, Academic Press, New York(1983).
\bibitem{shepley} L. C. Shepley and A. H. Taub, {\em Space-times containing perfect fluids and having a vanishing
conformal divergence}, Commun. Math. Phys(1967). vol.5,  pp-237–256.
\bibitem{mohd}Mohd Danish Siddiqi, {\em Conformal $\eta$-Ricci solitons in $\delta$- Lorentzian Trans Sasakian manifolds}, International Journal of Maps in Mathematics(2018), vol. 1, Isu. 1, pp- 15-34.
\bibitem{singh}Abhishek Singh, Shyam Kishor, {\em Some types of $\eta$-Ricci Solitons on Lorentzian para-Sasakian manifolds}, Facta Universitatis (NI\v{S})
\bibitem{soumendu} Soumendu Roy and Arindam~~Bhattacharyya, {\em Conformal Ricci solitons on 3-dimensional trans-Sasakian manifold}, Jordan Journal of Mathematics and Statistics (2020), Vol- 13(1), pp-89-109.
\bibitem{roy} Soumendu Roy, Santu Dey and Arindam Bhattacharyya, Shyamal Kumar Hui, {\em $*$-Conformal $\eta$-Ricci Soliton on Sasakian manifold}, 	arXiv:1909.01318v1 [math.DG] (2019).
\bibitem{roy2} Soumendu Roy, Santu Dey and Arindam Bhattacharyya, {\em Yamabe Solitons on $(LCS)_{n}$-manifolds}, arXiv:1909.06551v1 [math.DG] (2019).
 \bibitem{roy3} Soumendu Roy, Santu Dey and Arindam Bhattacharyya, {\em Some results on $\eta$-Yamabe Solitons in 3-dimensional trans-Sasakian manifold},  arXiv:2001.09271v2 [math.DG] (2020).
 \bibitem{stephani} H. Stephani, {\em General Relativity-An Introduction to the Theory of Gravitational Field} Cambridge University Press(1982.), Cambridge.
\bibitem{topping}Peter Topping, {\em Lecture on the Ricci Flow}, Cambridge University Press(2006).
\bibitem{yano} K. Yano, {\em On the torse-forming directions in Riemannian spaces}, Proc. Imp. Acad. Tokyo(1944), 20 , pp-340–345.
 \end{thebibliography}
\end{document}